\documentclass[10pt]{article}
\usepackage{amsmath}
\usepackage{amsfonts}
\usepackage{amssymb}
\usepackage{comment}
\newtheorem{theorem}{Theorem}[section]
\newtheorem{definition}[theorem]{Definition}
\newtheorem{proposition}[theorem]{Proposition}
\newtheorem{corollary}[theorem]{Corollary}
\newtheorem{lemma}[theorem]{Lemma}

\newcommand{\Proof}{{\it Proof. \/}}
\newcommand{\squareforqed}{\hbox{\rlap{$\sqcap$}$\sqcup$}}
\newcommand{\qed}{\ifmmode\squareforqed\else{\unskip\nobreak\hfil
\penalty50\hskip1em\null\nobreak\hfil\squareforqed
\parfillskip=0pt\finalhyphendemerits=0\endgraf}\fi}

\newcommand{\fp}{\qed\removelastskip\vskip\baselineskip\relax}
\DeclareMathOperator{\disc}{disc}

\DeclareMathOperator{\GL}{GL}

\newcommand{\psmm}[4]{\left(\begin{smallmatrix}{#1}&{#2}\\{#3}&{#4}\end{smallmatrix}\right)}
\newcommand{\Z}{\mathbb Z}
\newcommand{\Q}{\mathbb Q}
\newcommand{\R}{\mathbb R}

\newcommand{\al}{\alpha}
\newcommand{\be}{\beta}
\newcommand{\eps}{\varepsilon}

\newcommand{\ga}{\gamma}
\begin{document}
\pagestyle{plain}
\title{Polynomial Continued Fractions for Algebraic Numbers}
\author{Henri Cohen}

\maketitle

\begin{abstract}
  Our main result is that any real cubic algebraic number has a continued
  fraction expansion with polynomial coefficients. Some generalizations
  are mentioned.
\end{abstract}  
  
\section{Introduction}

Let $\al$ be a real algebraic number. If $\al$ is rational, it has a finite
simple continued fraction expansion, and if $\al$ is quadratic it has a
simple continued fraction expansion which is ultimately periodic.
The problem of studying the simple continued fraction expansion of higher
degree algebraic numbers, in particular cubic ones, has attracted a lot
of attention, and one can reasonably say that nothing is known. For instance
one does not know if the coefficients can be bounded, or if they are
always unbounded, etc...

Our goal is not to attack this notoriously difficult problem, but to
give a positive answer for certain algebraic numbers (in particular cubic ones)
if one generalizes the type of continued fraction that one
considers.

\begin{definition}\begin{enumerate}
  \item A \emph{generalized} continued fraction (abbreviated CF) is an
    expression of the type
  $$a(0)+b(0)/(a(1)+b(1)/(a(2)+b(2)/(a(3)+\cdots)))\;,$$
    where $a(n)$ and $b(n)$ are rational numbers, $b(n)\ne0$ for all $n$, and
    $a(n)\ne0$ except for a finite number of $n$.
\item Such a CF is said to be of \emph{polynomial type} if $a(n)$ and $b(n)$
  are polynomials for $n$ sufficiently large.\end{enumerate}
\end{definition}

{\bf Remarks.}\begin{enumerate}\item A \emph{simple} CF corresponds by
  definition to $b(n)=1$ for all $n$, $a(0)\in\Z$, and $a(n)\in\Z_{\ge1}$ for
  $n\ge1$.
\item By simplifying, it is clear that we can suppress all $a(n)=0$ and
  assume that $a(n)\ne0$ for all $n\ge1$.\end{enumerate}

\smallskip

{\bf Notation:} If $a(n)$ is a polynomial $A(n)$ in $n$ for $n>n_a$ and
$b(n)$ is a polynomial $B(n)$ in $n$ for $n>n_b$, we will denote the above CF
of polynomial type by the very practical notation
$$((a(0),a(1),\dotsc,a(n_a),A(n)),(b(0),b(1),\dotsc,b(n_b),B(n)))$$

\smallskip

Our main result is the following:

\begin{theorem}
  Any real cubic algebraic number has a continued fraction expansion of
  polynomial type.
\end{theorem}

Note that it has been known probably since the 19th century that
\emph{some} cubic irrationalities have such an expansion, for instance
the smallest positive root of $x^3-x+t=0$ when $t$ is sufficiently small,
but the point of the theorem is that it is true for \emph{all} cubic
irrationalities. The proof given here is quite simple, but involves a few
subtleties.

\section{Preliminary Lemmas}

Before giving the proof of this theorem, we need a series of lemmas.
  
\begin{definition} We will say that a real number $u$ is of \emph{continued
  fraction type} if it is equal to the limit of a convergent
  continued fraction $((a(n)),(b(n)))$ of polynomial type in the above sense,
  and we denote by ${\cal C}\subset\R$ the set of such constants.
\end{definition}

\begin{lemma}\label{lem01} If $\psmm{a}{b}{c}{d}\in \GL_2(\Q)$ and
  $u\in {\cal C}$ we have $(au+b)/(cu+d)\in{\cal C}$.\end{lemma}

\Proof Immediate: if $u=a(0)+b(0)/(a(1)+b(1)/(a(2)+\cdots))$, then
if $c\ne0$ we have
\begin{align*}\dfrac{au+b}{cu+d}&=\dfrac{a}{c}+\dfrac{-(ad-bc)/c^2}{u+d/c}\\
  &=(a/c)+(-(ad-bc)/c^2)/(a(0)+d/c+b(0)/(a(1)+b(1)/(a(2)+\cdots)))\;,\end{align*}
so is of the form $a'(0)+b'(0)/(a'(1)+\cdots)$ with $a'(0)=a/c$,
$a'(1)=a(0)+d/c$, $a'(n)=a(n-1)$ for $n\ge2$, $b'(0)=-(ad-bc/c^2)$, and
$b'(n)=b(n-1)$ for $n\ge1$. If $c=0$ the result is even more trivial since
$(au+b)/(cu+d)=(a/d)u+(b/d)$.\fp

\begin{lemma}\label{lem02} Let $K=\Q(u)$ be a cubic number field. Any
  element $v\in K\setminus\Q$ is of the form $v=(au+b)/(cu+d)$ for some
  $\psmm{a}{b}{c}{d}\in\GL_2(\Q)$.\end{lemma}

\Proof The four elements $\{uv,v,u,1\}$ belong to $K$ which is a $\Q$-vector
space of dimension $3$, so are linearly dependent, and since $v\notin\Q$
the corresponding $2\times2$ matrix is invertible, proving the lemma.\fp

It follows from these two lemmas that to prove our theorem it is sufficient
to prove that at least one non-rational element of a cubic number field is in
${\cal C}$.

We now set the following definition:

\begin{definition} Let $P$ be an irreducible polynomial of degree exactly
  equal to $3$, and write $P(x)=ax^3+bx^2+cx+d$. We define the \emph{ratio}
  of $P$ and denote by $r(P)$ the quantity
  $$r(P)=\dfrac{4(3ac-b^2)^3}{(27a^2d-9abc+2b^3)^2}$$
  By extension, if $u$ is an algebraic number of degree $3$ we will denote
  by $r(u)$ the ratio of the characteristic polynomial of $u$.
\end{definition}

For instance, if $P(x)=x^3+ax+b$ then we simply have $r(P)=4a^3/(27b^2)$,
and in particular if $P(x)=x^3-cx+c$ we have $r(P)=-4c/27$.

It is immediate to check that if $A$ and $B$ are rational with $A\ne0$
we have $r(Au+B)=r(u)$.

\smallskip

\begin{lemma}\label{lem1} If $\disc(P)>0$ we have $r(P)<-1$, and if
  $\disc(P)<0$ we have $r(P)>-1$.\end{lemma}

\Proof Indeed, we immediately check  that
$$r(P)+1=-\dfrac{27a^2\disc(P)}{(27a^2d-9abc+2b^3)^2}\;,$$
and the result follows.\fp

\begin{lemma}\label{lem2} If $K$ is a cubic number field then $|r(u)|$ is
  unbounded for $u\in K$.\end{lemma}

\Proof Let $K=\Q(\al)$ with $\al$ a root of a cubic polynomial which we may
assume to be $x^3+ax+b$, and set $u=\al^2+D\al+2a/3$. The characteristic
polynomial of $u$ is of the form $x^3+ex+f$ (i.e., again with no term in
$x^2$), with $e=aD^2+3bD-a^2/3$ and $f=bD^3-(2a^2/3)D^2-abD-(b^2+2a^3/27)$, so
$r(u)=4e^3/(27f^2)$. We compute that the resultant of these two polynomials in
$D$ is equal to $-\disc(\al)^3/729$, hence is nonzero, so it follows that if
we choose a rational value for $D$ close to a real root of the cubic polynomial
$f$, $|f|$ will be as small as we want while $|e|$ will be bounded from below,
so $r(u)$ is unbounded as desired.\fp

\begin{corollary}\label{cor2} Let $K$ be a cubic number field. There exists a
  defining equation for $K$ of the form $x^3-cx+c=0$ (so that $c>27/4$ if $K$
  is totally real and $c<27/4$ otherwise), and $|c|$ can be chosen as
  large as one wants.\end{corollary}

\Proof Clear: let $u$ be chosen so that $|r(u)|$ is large. If the
characteristic polynomial of $u$ is $x^3+a_2x^2+a_1x+a_0$, one checks
that $v=((9a_2^2-27a_1)/(27a_0+2a_2^3-9a_1a_2))(u+a_2/3)$ has a characteristic
polynomial of the form $x^3-cx+c$, so $|c|=(27/4)|r(v)|=(27/4)|r(u)|$
is as large as we like.\fp

\begin{lemma}\label{lem55} if $c>27/4$, the equation $x^3-cx+c=0$ has three
  real roots $\be_1$, $\be_2$, $\be_3$ which are such that
  $\be_1<-3<1<\be_2<3/2<\be_3$.
\end{lemma}

\Proof Trivial and left to the reader.\fp

\begin{lemma}\label{lem6} Let $c>27/4$, $c\ne 27/2$, and set
  $\phi_c(x)=3c(x-3)/((2c-27)x)$. Then the $\phi_c(\be_i)$ are the three
  roots of $x^3-Cx+C$ with $C=27c^2/(2c-27)^2$, and if we denote by
  $\ga_1<-3<1<\ga_2<3/2<\ga_3$ these three roots, we have:
  \begin{itemize}\item If $c<27/2$:
    $\phi_c(\be_1)=\ga_1$, $\phi_c(\be_2)=\ga_3$, and $\phi_c(\be_3)=\ga_2$.
  \item If $c>27/2$:
    $\phi_c(\be_1)=\ga_3$, $\phi_c(\be_2)=\ga_1$, and $\phi_c(\be_3)=\ga_2$.
  \end{itemize}
\end{lemma}
  
\Proof This is a simple direct verification and left to the reader. Note that
since in our application $P$ will be irreducible, we cannot have $c=27/2$
since in that case $x^3-cx+c$ is reducible.\fp

The main tool that we are going to use to prove our main theorem is the
following crucial lemma, which is a direct consequence of the Lagrange
inversion formula:

\begin{lemma}\label{lem4}\begin{enumerate}
  \item Assume that $|c|>27/4$. The series
    $$S(c)=\sum_{n\ge0}\dfrac{1}{2n+1}\binom{3n}{n}c^{-n}$$
    converges to a root $\be=S(c)$ of $x^3-cx+c=0$.
  \item More precisely, if $c<-27/4$, $\be$ is its unique real root, while
    if $c>27/4$, $\be=\be_2$ is its unique root such that $1<\be<3/2$.
  \item More generally, let $P=x^3+ax+b$. If $|r(P)|=|4a^3/(27b^2)|>1$, one
    of the roots $\al$ of $P=0$ is given by
    $$\al=(-b/a)\sum_{n\ge0}\dfrac{1}{2n+1}\binom{3n}{n}(-b^2/a^3)^n\;,$$
    and $\al$ is the unique root of $P$ if $4a^3/(27b^2)>1$, and otherwise
    $\al$ is the unique root of $P$ such that $|b/a|<|\al|<(3/2)|b/a|$.
  \end{enumerate}
\end{lemma}

\Proof (1). This is well-known and is a simple application of Lagrange's
inversion formula. Note that the simplicitly of the result comes from the fact
that if $f(x)=x^3-cx+c$, then $(x/(f(x)-f(0)))^n=(x^2-c)^{-n}$ has a simple
power series expansion around $x=0$.

\smallskip

(2). If $c<-27/4$, $P$ has a unique real root which therefore must be
equal to $\be$, so assume that $c>27/4$. I claim that $0<\be-1<(27/(8c))$,
and since $c>27/4$ we thus have $1<\be<3/2$, which by Lemma \ref{lem55}
is the unique root in that interval, proving (2). To prove my claim,
by (1) we have $\be^3-c\be+c=0$, so $c(\be-1)=\be^3$. On the other hand, in
view of the representation as a series it is clear that $\be>1$ and since
$c>27/4$ we have $\be<\sum_{n\ge0}(\binom{3n}{n}/(2n+1))(27/4)^{-n}$. It is
immediate to check that this series is still convergent, and by continuity is
a root of $x^3-(27/4)x+(27/4)=0$ which are $-3$ and $3/2$, so its sum is equal
to $3/2$. Thus $c(\be-1)=\be^3<(3/2)^3$, proving our claim, hence (2).

\smallskip

(3). Follows immediately from (1) by setting $\be=-(a/b)\al$.\fp

\begin{corollary}\label{cor3} The elements $\al$ and $\be$ given by Lemma
  \ref{lem4} are of continued fraction type, in other words belong to
  ${\cal C}$.
\end{corollary}

\Proof This immediately follows from Euler's transformation of a series
into a continued fraction. To be completely explicit, using the short-hand
notation introduced above we have
$$\be=((1,6c,(4c+27)n^2+(2c-27)n+6),(6,-6cn(2n+1)(3n+1)(3n+2)))\;,$$ in
other words
$$\be=1+6/(6c-360c/(20c+60-3360c/(42c+168-\cdots)))\;,$$
with $a(0)=1$, $a(1)=6c$, $a(n)=(4c+27)n^2+(2c-27)n+6$ for $n\ge2$, $b(0)=6$,
and $b(n)=-6cn(2n+1)(3n+1)(3n+2)$ for $n\ge1$, and similarly for $\al$.\fp

\section{Proof of the Theorem}

We can now restate and prove the following:

\begin{theorem}\label{thm1}\begin{enumerate}\item
    If $u$ is a real algebraic number of degree $3$ then $u\in{\cal C}$.
  \item More precisely, there exists a constant $E$ with $|E|>1$ and a CF of
    polynomial type with convergents $p(n)/q(n)$ converging to $u$ such that
    $$u-\dfrac{p(n)}{q(n)}\sim\dfrac{C}{E^nn^{3/2}}$$
    for a suitable constant $C\ne0$, and $|E|$ can be taken arbitrarily large.
    \end{enumerate}
  \end{theorem}

\Proof (1). Let $K=\Q(u)$ be the cubic number field generated by $u$. By
Corollary \ref{cor2} we can choose a defining equation for $K$ of the
form $x^3-cx+c=0$ with $|c|$ as large as we like, in particular for now
$|c|>27/4$, so $K=\Q(v)$ with $v$ a real root of $P(x)=x^3-cx+c=0$.
If $K$ is not totally real, $v$ is the unique real root of $P$, so
by Lemma \ref{lem4} is the sum of the convergent hypergeometric series,
so belongs to ${\cal C}$ by Corollary \ref{cor3}. Since by Lemma \ref{lem02}
$u$ is the image of $v$ by $\GL_2(\Q)$ and ${\cal C}$ is stable under this
action, it follows that $u\in{\cal C}$.

Assume now that $K$ is totally real, so that $c>27/4$. Since we can choose
$|c|$ as large as we desire, we now choose $c>27/2$. The polynomial
$x^3-cx+c=0$ now has three real roots $\be_i$ as in Lemma \ref{lem55}, and
by Lemma \ref{lem4} the sum of the hypergeometric series is equal to
$\be_2$. Thus if $1<v<3/2$ we have $v=\be_2$ and we are done as in the
non-totally real case. If $v>3/2$, we have $v=\be_3$ hence by Lemma \ref{lem6}
$\phi_c(v)=\ga_2$ is a root of some polynomial $x^3-Cx+C$, so is the sum of the
corresponding hypergeometric series, hence $\phi_c(v)\in{\cal C}$, and since
$\phi_c$ is a $\GL_2(\Q)$ action and ${\cal C}$ is stable, it follows once
again that $v\in{\cal C}$. Finally, assume that $v<-3$, so that $v=\be_1$.
Since $c>27/2$, by Lemma \ref{lem6} once again, we have $\phi_c(v)=\ga_3$, the
unique root of $x^3-Cx+C=0$ for $C=27c^2/(2c-27)^2$ with $\ga_3>3/2$.
We easily check that $C>c>27/2$, so again applying Lemma \ref{lem6} we deduce
that $\phi_C(\ga_3)=\delta_2$, where $\delta_2$ is the middle root of still
another equation of the type $x^3-C_1x+C_1=0$, so $\delta_2\in{\cal C}$, and
thus $v\in{\cal C}$ also in this case, proving (1).

\smallskip

(2). By Euler's transformation of series into CFs and estimating the remainder
of a convergent hypergeometric series, it is immediate to show that with the
notation of Lemma \ref{lem4}, we have
$$S(c)-\dfrac{p(n)}{q(n)}\sim\dfrac{C}{(4c/27)^nn^{3/2}}$$
for a suitable constant $C\ne0$ depending on $c$. Applying an element of
$\GL_2(\Q)$ to $S(c)$ only changes the constant $C$ but not the main
asymptotics, proving (2) since by Corollary \ref{cor2} for a given $u$
we can choose $|c|$ as large as we like.\fp

\section{Examples}

To simplify notation, as in Lemma \ref{lem4}, for $|c|>27/4$ we will set
$$S(c)=\sum_{n\ge0}\dfrac{1}{2n+1}\binom{3n}{n}c^{-n}\;,$$
which by Euler has the CF expansion given in the proof of Corollary \ref{cor3},
which we write in abbreviated form as
$$S(c)=((1,6c,(4c+27)n^2-(27c-2)n+6c),(6c,-6cn(2n+1)(3n+1)(3n+2)))$$
Every cubic irrationality $\be$ will be of the form $\be=(AS(c)+B)/(CS(c)+D)$
for a suitable $c$ with $|c|>27/4$ and $\psmm{A}{B}{C}{D}\in\GL_2(\Q)$.
If $C=0$ we have $\be=(A/D)S(c)+B/D$, giving the CF
\begin{align*}\be&=(((A+B)/D,6c,(4c+27)n^2-(27c-2)n+6c),\\&\phantom{=}(6(A/D)c,-6cn(2n+1)(3n+1)(3n+2)))\;.\end{align*}
If $C\ne0$, as explained above we write $\be=A/C+(-(AD-BC)/C)/(S(c)+D/C)$,
giving the CF
\begin{align*}\be&=((A/C,1+D/C,6c,(4c+27)n^2-(6c+81)n+2c+60),\\&\phantom{=}((-(AD-BC)/C),6c,-6c(n-1)(2n-1)(3n-1)(3n-2)))\end{align*}
Thus, the CF expansion of $\be$ is entirely determined by $c$ and
$\psmm{A}{B}{C}{D}$, which is therefore the only data that we will give
in the form $\be=(c,\psmm{A}{B}{C}{D})$.

\smallskip

$\bullet$ {\bf Discriminant $-23$}.

\smallskip

An equation for the cubic field of discriminant $-23$ is $x^3+x^2+2x+1=0$,
whose unique real root is $\be=-0.56984029\cdots$. In this case, we already
have $r(\be)=500/121>1$, and we deduce that
$$\be=\left(-\dfrac{3375}{121},\begin{pmatrix}-11/45&-1/3\\0&1\end{pmatrix}\right)\;.$$

\smallskip

$\bullet$ {\bf Discriminant $-31$}.

\smallskip

An equation for the cubic field of discriminant $-23$ is $x^3+x+1=0$,
whose unique real root is $\be=-0.68232780\cdots$. In this case, we have
$r(\be)=4/27<1$ so we need to work some more. Setting $\al=-\be$ and
$u=\al^2+D\al+2/3$ as in the proof of Lemma \ref{lem2}, we find that for
$D=-1$ we have $r(u)=5324/2209>1$, giving a hypergeometric representation
of $u$ as in Lemma \ref{lem4} with $a=11/3$ and $b=-47/27$, hence a CF
using Euler's transformation. We must now express $\al$ as a $\GL_2(\Q)$
image of $u$, and we find that $\al=(3u-5)/(3u+4)$, and we deduce that
$$\be=\left(-\dfrac{35937}{2209},\begin{pmatrix}47&-165\\47&132\end{pmatrix}\right)\;.$$

\smallskip

$\bullet$ {\bf Discriminant 49}.

\smallskip

We now have a totally real \emph{cyclic} cubic field, so the situation is
simpler since given one root of a defining polynomial, the other roots are
obtained by a M\"obius transformation, so have the same CF expansion apart
from the inital terms.
An equation for this field is $x^3-x^2-2x+1=0$, which has the
three roots $\be_1=-1.246979\cdots$, $\be_2=0.445041\cdots$, and
$\be_3=1.8019377\cdots$. In this particularly simple case we can apply
Lemma \ref{lem4} directly, and obtain the following:
$$\be_1=\left(189,\begin{pmatrix}1&-6\\1&3\end{pmatrix}\right)\;,\quad
  \be_2=\left(189,\begin{pmatrix}1/9&1/3\\0&1\end{pmatrix}\right)\;,\quad
      \be_3=\left(189,\begin{pmatrix}0&-9\\1&-6\end{pmatrix}\right)\;.$$
        
\smallskip

$\bullet$ {\bf Discriminant 148}.

\smallskip

A defining equation for this field is $x^3+x^2-3x-1=0$, which has the three
roots $\be_1=-2.1700864\cdots$, $\be_2=-0.3111078\cdots$, and
$\be_3=1.4811943\cdots$. We compute that $r(\be_i)=-1000<-1$, so we can
apply Lemma \ref{lem4} directly, which will give the middle root $\be_2$,
and we find that
$$\be_2=\left(6750,\begin{pmatrix}1/45&-1/3\\0&1\end{pmatrix}\right)\;.$$
We now apply the proof of the theorem which tells us that a root of
$x^3-Cx+C=0$ with $C=27c^2/(2c-27)^2$ will be in the number field generated
by $\be_3$. We find that
$$\be_3=\left(\dfrac{1687500}{249001},\begin{pmatrix}-499&600\\1497&-2250\end{pmatrix}\right)\;.$$
We iterate this with $C_1=27C^2/(2C-27)^2$, and we find that
$$\be_1=\left(\dfrac{105468750000}{15376248001},\begin{pmatrix}-496004&1121250\\1860015&-2801250\end{pmatrix}\right)\;.$$

Note that it is of course quite possible that there are expansions with
simpler coefficients.

\smallskip

$\bullet$ {\bf The case $\be=\root3\of2$}.

\smallskip

The minimal polynomial is of course $x^3-2$, which we must change so as to
have $a\ne0$, and as in the proof of Lemma \ref{lem2} we replace $\be$ by
$\be^2+\be$, but in fact even better by $\al=-(\be^2+\be)$ whose characteristic
polynomial is $x^3-6x+6$, already of the desired shape. Looking as before for
$D$ such that $\ga=\al^2+D\al-4$ has $r(\ga)>1$, we find that we can take
$D=3/2$, and we immediately find the representation
$$\root3\of2=\left(-54,\begin{pmatrix}-1&6\\1&3\end{pmatrix}\right)\;.$$
  
\section{An Easy Result}

\begin{proposition} If $u$ and $a$ are rational numbers and $u>0$, then
  $u^a$ is the sum of a hypergeometric series. In particular $u^a\in{\cal C}$.
\end{proposition}

\Proof Write $a=m/d$ with $d>0$, and consider $\al=u^{-1/d}$.
If $\al\in\Q$ then $u^a=\al^{-m}\in\Q$ so the result is clear. Otherwise,
let $f$ be a rational approximation to $\al$, and write $\eps=f-\al$.
We have
$$uf^d-1=u(\al+\eps)^d-1=u(\al^d+\eps d\al^{d-1}+O(\eps^2))-1=\eps \al d+O(\eps^2)\;,$$
so for $\eps$ sufficiently small (essentially $\eps<1/(2\al d)$), we may
assume that $0<z=uf^d-1<1$. Since $(1+z)^a=u^af^m$, we thus have the
convergent hypergeometric representation
$u^a=f^{-m}\sum_{n\ge0}\binom{a}{n}z^n$.\fp

For example, the following formulas give hypergeometric representations of
$\root3\of2$:

\begin{align*}\root3\of2&=(1-1/2)^{-1/3}=(4/3)(1+5/27)^{-1/3}\\
  &=(5/4)(1+3/125)^{1/3}=(5/4)(1-3/128)^{-1/3}=\cdots\end{align*}
For instance, the very first formula gives the polynomial continued fraction
representation
$$\root3\of2=((1,6,9n-2),(1,-6n(3n+1)))\;.$$
All these CFs are of course of a much simpler form than those obtained
as above from the general results on cubic irrationalities.

\section{Generalizations}

We were able to treat the cubic case because of two (related) facts: first
and most importantly, Lagrange inversion gave a hypergeometric series thanks
to the reduction of the defining equation to the trinomial $x^3+ax+b$,
and second because two irrational elements of a cubic field are related by
a $\GL_2(\Q)$ transformation.

A generalization of Lemma \ref{lem4} to general degrees is as follows,
again a direct consequence of Lagrange inversion:

\begin{proposition} Let $P=x^k+ax+b$. If $|(k-1)^{k-1}a^k/(k^kb^{k-1})|>1$,
  one of the roots $\al$ of $P=0$ is given by
  $$\al=(-b/a)\sum_{n\ge0}\dfrac{1}{(k-1)n+1}\binom{kn}{n}((-1)^kb^{k-1}/a^k)^n\;.$$
\end{proposition}

Using this and methods analogous to those in the cubic case, one can prove
that many algebraic numbers of degree greater or equal to 4 are in ${\cal C}$,
but I do not know a general statement.

\end{document}